\documentclass {article}

 \usepackage [english] {babel}

 \usepackage {amssymb, amsmath}
 \usepackage {xypic}

 \def \:{\colon}
 \def \le {\leqslant}
 \def \ge {\geqslant}
 \def \[{[\![}
 \def \]{]\!]}
 \def \Z{\boldsymbol{\mathsf Z}}
 \def \R{\boldsymbol{\mathsf R}}
 \def \d{\partial}

 \def \break{\par \addvspace {\medskipamount} \par}

 \newenvironment {claim} [1]
     {\break \noindent {\bf #1} \it}
     {\break}
 \newenvironment {demo} [1]
     {\break \noindent {\it #1} }
     {\quad $\square$ \break}

 \DeclareMathOperator {\id} {id}
 \DeclareMathOperator {\incl} {in}
 \DeclareMathOperator {\proj} {pr}
 \DeclareMathOperator {\sk} {sk}
 \DeclareMathOperator {\cosk} {cosk}
 \DeclareMathOperator {\dist} {dist}
 \DeclareMathOperator {\diam} {diam}

\begin {document} \frenchspacing

 \title {An iterated sum formula \\
         for a spheroid's homotopy class modulo 2--torsion}

 \author {S. S. Podkorytov}

 \date {}

 \maketitle

 \begin {abstract} \noindent
 Let $X$ be a simply connected pointed space with finitely
 generated homotopy groups.
 Let $\Pi_n(X)$ denote the set of all continuous maps
 $a\:I^n\to X$ taking $\d I^n$ to the basepoint.
 For $a\in\Pi_n(X)$,
 let $[a]\in\pi_n(X)$ be its homotopy class.
 For an open set $E\subset I^n$,
 let $\Pi(E,X)$ be the set of all continuous maps $a\:E\to X$
 taking $E\cap\d I^n$ to the basepoint.
 For a cover $\Gamma$ of $I^n$,
 let $\Gamma(r)$ be the set of all unions of at most $r$
 elements of $\Gamma$.
 Put $r=(n-1)!$.
 We prove that
 for any finite open cover $\Gamma$ of $I^n$
 there exist maps $f_E\:\Pi(E,X)\to\pi_n(X)\otimes\Z[1/2]$,
 $E\in\Gamma(r)$, such that
 $$
 [a]\otimes1=\sum_{E\in\Gamma(r)}f_E(a|_E)
 $$
 for all $a\in\Pi_n(X)$.
 \end {abstract}


 \section * {1. Introduction}

 Let $X$ be a (pointed) space.
 An \textit {$n$--spheroid} in $X$ is a (continuous) map
 $a\:I^n\to X$ taking $\d I^n$ to the basepoint.
 Let $\Pi_n(X)$ denote the set of all $n$--spheroids in $X$.
 For $a\in\Pi_n(X)$,
 let $[a]\in\pi_n(X)$ be its homotopy class.
 For an open set $E\subset I^n$,
 let $\Pi(E,X)$ be the set of all maps $a\:E\to X$ taking
 $E\cap\d I^n$ to the basepoint.
 For a cover $\Gamma$ of $I^n$,
 let $\Gamma(r)$ be the set of all unions of at most $r$
 elements of $\Gamma$.

 Let $L$ be an abelian group.
 Consider a functional (i.~e. a function)
 $f\:\Pi_n(X)\to L$.
 We define the \textit {degree} of $f$ (denoted $\deg f$) to be
 the infimum of $r$ such that
 for any finite open cover $\Gamma$ of $I^n$
 there exist functionals $f_E\:\Pi(E,X)\to L$, $E\in\Gamma(r)$,
 such that
 $$
 f(a)=\sum_{E\in\Gamma(r)}f_E(a|_E)
 $$
 for all $a\in\Pi_n(X)$.
 We are interested in functionals of finite degree.

 This definition is motivated by the notion of order-restricted
 perceptron \cite {MP}.
 For relation to Vassiliev knot invariants, see \cite {Z}.
 An example of a functional $f\:\Pi_n(X)\to\R$ with $\deg\le r$
 is given by
 $$
 f(a)=\int_{u_1,\ldots,u_r\in I^n}
          p(u_1, a(u_1),\ldots,u_r,a(u_r))\,
          du_1\ldots du_r,
 $$
 where $p\:(I^n\times X)^r\to\R$ is any decent (say, bounded)
 measurable function.

 We prove the following result.

 \begin {claim} {1.1. Theorem.}
 Suppose $X$ is simply connected and has finitely generated
 homotopy groups.
 Define $q\:\Pi_n(X)\to\pi_n(X)\otimes\Z[1/2]$ by
 $q(a)=[a]\otimes1$.
 Then $\deg q\le(n-1)!$.
 \end {claim}

 This may be viewed as a finite sum version of Chen's iterated
 integrals.
 For a related result, see \cite {me}.
 Possibly, the $\Z[1/2]$ factor and the factorial sign may be
 removed.
 Claim~1.2 implies that $\deg q\ge n-1$ for $X=S^2\vee S^2$.

 \begin {claim} {1.2. Claim.}
 Let $f\:\Pi_n(X)\to L$ be a homotopy invariant functional with
 $\deg f\le r$.
 Let $b_0,b\in\Pi_n(X)$ be spheroids with
 $$
 [b_0]=0, \quad
 [b]=\[\ldots\[\[v_0,v_1\],v_2\],\ldots,v_r\],
 $$
 where $v_s\in\pi_{k_s}(X)$
 ($k_0+\ldots+k_r=n+r$;
 $\[{\cdot},{\cdot}\]$ denotes the Whitehead product).
 Then $f(b)=f(b_0)$.
 \end {claim}

 (A functional $f\:\Pi_n(X)\to L$ is \textit {homotopy
 invariant} if $f(a)$ depends only on $[a]$.)

 \begin {demo} {Proof \rm (cf. \cite [Theorem~3.2] {MP}).}
 Put
 $$
 T=S^{k_0}\vee\ldots\vee S^{k_r}.
 $$
 Let $g\:T\to X$ be a map whose restriction to the $s$th wedge
 summand represents $v_s$.
 Let $e_s\in\pi_{k_s}(T)$ be the elements represented by the
 canonical embeddings $S^{k_s}\to T$.
 Let $z\in\Pi_n(T)$ be a spheroid with
 $[z]=\[\ldots\[\[e_0,e_1\],e_2\],\ldots,e_r\]$.
 We have $[g\circ z]=[b]$.
 For every $k$,
 choose maps $p^k_i\:S^k\to S^k$, $i=0,1$, such that
 $p^k_i$ has degree $i$ and
 $p^k_0|_{U^k}=p^k_1|_{U^k}$
 for some open neighbourhood $U^k$ of the basepoint.
 Put
 $$
 p_{i_0,\ldots,i_r}=
 p^{k_0}_{i_0}\vee\ldots\vee p^{k_r}_{i_r}\:T\to T,
 \qquad i_0,\ldots,i_r=0,1.
 $$
 We have
 $$
 [p_{i_0,\ldots,i_r}\circ z]=i_0\ldots i_r[z]
 $$
 in $\pi_n(T)$
 and thus
 $$
 [g\circ p_{i_0,\ldots,i_r}\circ z]=i_0\ldots i_r[b]
 $$
 in $\pi_n(X)$.
 The sets
 $$
 V_s=U^{k_0}\vee\ldots\vee S^{k_s}\vee\ldots\vee U^{k_r},
 \qquad s=0,\ldots,r,
 $$
 form an open cover of $T$.
 Put $\Gamma=\{\,z^{-1}(V_s)\mid s=0,\ldots,r\,\}$.
 Since $\deg f\le r$,
 there are functionals $f_E\:\Pi(E,X)\to L$, $E\in\Gamma(r)$,
 such that
 $$
 f(a)=\sum_{E\in\Gamma(r)}f_E(a|_E)
 $$
 for all $a\in\Pi_n(X)$.
 We have
 \begin {multline*}
 (-1)^r(f(b_0)-f(b))=
 \sum_{i_0,\ldots,i_r=0,1}
 (-1)^{i_0+\ldots+i_r}
 f(g\circ p_{i_0,\ldots,i_r}\circ z)= \\
 =\sum_{E\in\Gamma(r)}
 \sum_{i_0,\ldots,i_r=0,1}
 (-1)^{i_0+\ldots+i_r}
 f_E(g\circ p_{i_0,\ldots,i_r}\circ z|_E).
 \end {multline*}
 Take $E\in\Gamma(r)$.
 We have
 $$
 z(E)\subset
 S^{k_0}\vee\ldots\vee U^{k_s}\vee\ldots\vee S^{k_r}
 $$
 for some $s$.
 Thus $p_{i_0,\ldots,i_r}\circ z|_E$ does not depend on $i_s$.
 Thus the inner sum equals zero.
 \end {demo}

 \paragraph * {Conventions and notation.}
 \textit {Maps} are continuous,
 unlike \textit {functions}.
 A \textit {space} is a pointed space.
 A \textit {subspace} contains the basepoint.
 Maps between spaces are basepoint preserving.
 This applies also to homotopies etc.
 A \textit {cell space} is a pointed CW complex.
 $\sk_n X$ denotes the $n$-skeleton of a cell space $X$.

 The homotopy relation is denoted $\sim$;
 $\sim_A$ is used for homotopy rel $A$.
 For a pair $(X,A)$,
 $\incl_{(X,A)}\:A\to X$ and $\proj_{(X,A)}\:X\to X/A$ are the
 inclusion and the projection.
 The subscript of $\incl$ and $\proj$ is often omitted.


 \section * {2. Making a loop space simply connected}

 The aim of this section is to prove Corollary~2.11.
 First fix some notation.

 \paragraph * {Homotopy fibres and cofibres.}
 Let
 $X$ and $Y$ be spaces,
 $g\:X\to Y$ be a map.
 Then we have the homotopy fibre sequence
 $$
 \xymatrix {
 F(g)
 \ar[r]^-{p(g)} &
 X
 \ar[r]^-g &
 Y,
 }
 $$
 where
 $F(g)=\{\,(x,v)\in X\times PY\mid g(x)=v(1)\,\}$ and
 $p(g)$ is the fibraton defined by $p(g)(x,v)=x$.
 We have the homotopy cofibre sequence
 $$
 \xymatrix {
 X
 \ar[r]^-g &
 Y
 \ar[r]^-{i(g)} &
 C(g),
 }
 $$
 where
 $C(g)$ is the unreduced cone of $g$ and
 $i(g)$ is the canonical embedding.


 \subsection * {2.A. A Moore space}

 Fix $d>0$.
 Let $M_d$ be the space obtained from $S^1$ by attaching a
 2-cell via a map $S^1\to S^1$ of degree $d$.
 Our aim here is to prove Corollary~2.5.

 \begin {claim} {2.1. Lemma.}
 Let
 $$
 \xymatrix {
 F
 \ar[r]^-i &
 E
 \ar[r]^-p &
 B
 }
 $$
 be a fibre sequence.
 Suppose
 $B$ and $E$ are path-connected,
 $p$ induces an isomorphism on $\pi_1$, and
 the canonical action of $\pi_1(B)$ on $\pi_2(B)$ is trivial.
 Then the canonical action of $\pi_1(B)$ on $H_1(F)$ is
 trivial.
 \end {claim}

 \begin {demo} {Proof.}
 The group $\pi_1(E)$ acts canonically on $\pi_1(F)$
 (see \cite [5.1.7.3] {FR});
 it also acts on $\pi_2(B)$ and $H_1(F)$ through
 $p_*\:\pi_1(E)\to\pi_1(B)$.
 The boundary homomorphism $\d$ and the Hurewicz homomorphism
 $h$
 $$
 \xymatrix {
 \pi_2(B)
 \ar[r]^-\d &
 \pi_1(F)
 \ar[r]^-h &
 H_1(F)
 }
 $$
 respect these actions
 (regarding $\d$, see \cite [5.1.8.4] {FR}).
 Since $\d$, $h$, and $p_*$ are epimorphisms,
 $\pi_1(B)$ acts trivially on $H_1(F)$.
 \end {demo}

 \begin {claim} {2.2. Claim.}
 $H_2(\Omega\Sigma M_d)\cong\Z_d$.
 \end {claim}

 \begin {demo} {}
 Easily seen from the homology spectral sequence of the path
 fibration of $\Sigma M_d$.
 \end {demo}

 \begin {claim} {2.3. Claim.}
 If $d$ is odd, $\pi_3(\Sigma M_d)\cong\Z_d$.
 \end {claim}

 \begin {demo} {Proof.}
 There is a fibration $U\to\Sigma M_d$ with fibre of homotopy
 type $K(\Z_d,1)$ and $U$ 2--connected.
 We have $\pi_3(\Sigma M_d)\cong H_3(U)$.
 Using the cohomology spectral sequence,
 we get\footnote [1]
 {This does not work for $d$ even.
 This results in the $\Z[1/2]$ factor in Theorem~1.1.}
 $H_3(U)\cong\Z_d$.
 \end {demo}

 \begin {claim} {2.4. Claim.}
 Suppose $d$ is odd.
 Then the canonical embedding $j\:M_d\to\Omega\Sigma M_d$
 induces zero homomorphism on $\pi_2$.
 \end {claim}

 \begin {demo} {Proof.}
 Consider the homotopy fibre sequence
 $$
 \xymatrix {
 F(j)
 \ar[r]^-{p(j)} &
 M_d
 \ar[r]^-j &
 \Omega\Sigma M_d.
 }
 $$
 By Lemma~2.1,
 the action of $\pi_1(\Omega\Sigma M_d)$ on $H_1(F(j))$ is
 trivial.
 Using the homology spectral sequence and Claim~2.2,
 we get $H_1(F(j))\cong\Z_d$.
 The boundary homomorphism
 $\d\:\pi_2(\Omega\Sigma M_d)\to\pi_1(F(j))$ is an epimorhism.
 Thus $\pi_1(F(j))$ is abelian and thus isomorphic to $\Z_d$.
 By Claim~2.3,
 $\pi_2(\Omega\Sigma M_d)\cong\Z_d$.
 Thus $\d$ is an isomorphism.
 Thus $j_*\:\pi_2(M_d)\to\pi_2(\Omega\Sigma M_d)$ is zero.
 \end {demo}

 \begin {claim} {2.5. Corollary.}
 Suppose $d$ is odd.
 Let
 $X$ be a space,
 $f\:M_d\to\Omega X$ be a map.
 Then $f$ induces zero homomorphism on $\pi_2$.
 \end {claim}

 \begin {demo} {Proof.}
 Let $g\:\Sigma M_d\to X$ be the map adjoint to $f$.
 Then $f=\Omega g\circ j$,
 where $j$ is as in Claim~2.4.
 \end {demo}


 \subsection * {2.B. Two technical lemmas}

 \begin {claim} {2.6. Lemma.}
 Let
 $X$ be a space, 
 $A\subset X$ be a closed subspace such that
 $(X,A)$ is a Borsuk pair.
 Suppose $X$ and $A$ are homotopy equivalent to cell spaces of
 dimension at most $n$ and $n-1$, respectively.
 Let
 $Y$ be a cell space
 and $f\:X\to Y$ be a map with $f(A)\subset\sk_nY$.
 Then $f$ is rel $A$ homotopic to a map $h\:X\to Y$ with
 $h(X)\subset\sk_nY$.
 \end {claim}

 This is used only for $n=3$.

 \begin {demo} {Proof.}
 The map $f$ is homotopic to a map $g\:X\to Y$ with
 $g(X)\subset\sk_nY$.
 Let $k\:X\times[0,1]\to Y$ be the corresponding homotopy:
 $$
 k(x,0)=f(x), \quad
 k(x,1)=g(x),
 \qquad x\in X.
 $$
 The map $k|_{A\times[0,1]}$ is rel $A\times\{0,1\}$ homotopic
 to a map $q\:A\times[0,1]\to Y$ with image in $\sk_nY$.
 Since $(X,A)$ is a Borsuk pair,
 there is a map $l\:X\times[1,2]\to Y$ with image in $\sk_nY$
 such that
 $l(x,1)=g(x)$, $x\in X$, and
 $l(x,1+t)=q(x,1-t)$, $x\in A$, $t\in[0,1]$.
 Let $m\:X\times[0,2]\to Y$ be the map with
 $m|_{X\times[0,1]}=k$, $m|_{X\times[1,2]}=l$.
 The map $m|_{A\times[0,2]}$ is rel $A\times\{0,2\}$ homotopic
 to the map $e\:A\times[0,2]\to Y$, $e(x,t)=f(x)$.
 By the Str{\o}m theorem \cite [Lecture~I, Proposition~2] {P},
 $$
 (X\times[0,2],X\times\{0,2\}\cup A\times[0,2])
 $$
 is a Borsuk pair.
 Thus there is a map $\tilde m\:X\times[0,2]\to Y$ with
 $\tilde m|_{X\times\{0,2\}}=m|_{X\times\{0,2\}}$,
 $\tilde m|_{A\times[0,2]}=e$.
 Put $h(x)=\tilde m(x,2)$.
 Then $\tilde m$ is the desired homotopy.
 \end {demo}

 \begin {claim} {2.7. Lemma.}
 Let
 $V$ be a Hausdorff space,
 $T\subset V$ be a compact subspace.
 Let
 $W$ be a space,
 $a\:V\to W$, $b\:V/T\to W$, $c\:T\to W$ be maps such that
 $b\circ\proj_{(V,T)}=a$, $c=a|_T$
 (so $c$ is constant).
 Let $j\:F(c)\to F(a)$ and $l\:F(a)\to F(b)$ be the maps
 induced by $\incl_{(V,T)}$ and $\proj_{(V,T)}$, respectively.
 Let
 $Y$ be a space,
 $X\subset Y$ be a subspace,
 $s\:F(a)\to X$ and $t\:F(b)\to Y$ be maps such that
 $\incl_{(Y,X)}\circ s$ and $t\circ l$ are rel $j(F(c))$
 homotopic.
 Then there exists a unique map $r\:F(b)\to X$ such that
 $r\circ l=s$.
 The maps $\incl_{(Y,X)}\circ r$ and $t$ are homotopic rel
 $p(b)^{-1}(q_0)$,
 where $q_0\in V/T$ is the basepoint.
 \end {claim}

 $$
 
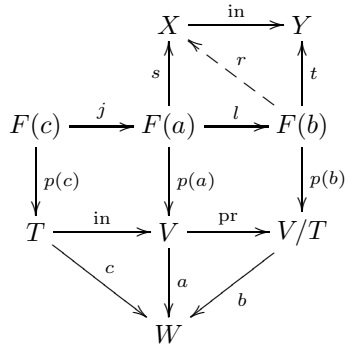
\xymatrix {
 &
 X
 \ar[r]^-\incl &
 Y  \\
 F(c)
 \ar[r]^-j
 \ar[d]^-{p(c)} &
 F(a)
 \ar[u]^-s
 \ar[r]^-l
 \ar[d]^-{p(a)} &
 F(b)
 \ar@{-->}[ul]_-r
 \ar[u]_-t
 \ar[d]^-{p(b)} \\
 T
 \ar[r]^-\incl
 \ar[dr]^-c &
 V
 \ar[r]^-\proj
 \ar[d]^-a &
 V/T
 \ar[dl]^-b \\
 &
 W &
 }
 $$

 \begin {demo} {}
 This is because $l$ is a quotient map.
 \end {demo}


 \subsection * {2.C. A $K(G,1)$ space
                and some its subquotients}

 Let
 $G$ be a finitely generated abelian group and
 $$
 G=G_1\oplus\ldots\oplus G_r
 $$
 be its decomposition into cyclic summands.
 For each $s=1,\ldots,r$,
 take a cell space $V_s$ of homotopy type $K(G_s,1)$ such that
 $\sk_2V_s$ is either $S^1$ or $M_d$ (see 2.A) for proper $d$.
 Put
 $T_s=\sk_1V_s$ ($=S^1$),
 $U_s=\sk_3V_s$,
 $$
 V=V_1\times\ldots\times V_r, \quad
 T=T_1\times\ldots\times T_r, \quad
 U=U_1\times\ldots\times U_r.
 $$
 We have $T\subset U\subset V$.
 $V$ is a $K(G,1)$ space.

 \begin {claim} {2.8. Claim.}
 Let
 $X$ be a simply connected space with $\pi_2(X)\cong G$,
 $Q$ be a space of weak homotopy type $K(G,2)$,
 $b\:X\to Q$ be a map inducing an isomorphism on $\pi_2$, and
 $m\:V\to\Omega Q$ be a weak homotopy equivalence.
 Suppose $G$ has no 2--torsion.
 Then there exists a map $f\:U\to\Omega X$ such that
 $\Omega b\circ f\sim m|_U$.
 \end {claim}

 $$
 \xymatrix {
 U
 \ar@{-->}[d]_-f
 \ar[r]^-\incl &
 V
 \ar[d]^-m \\
 \Omega X
 \ar[r]^-{\Omega b} &
 \Omega Q
 }
 $$

 \begin {demo} {Proof.}
 We naturally have $U_s\subset U$, $s=1,\ldots,r$.
 For every $s$,
 there is a map $f'_s\:\sk_2U_s\to\Omega X$ such that
 $\Omega b\circ f'_s$ and $m|_{\sk_2U_s}$ induce the same
 homomorphism on $\pi_1$.
 By Corollary~2.5,
 $f'_s$ induces zero homomorphism on $\pi_2$.
 Thus it extends to a map $f_s\:U_s\to\Omega X$.
 Since $\Omega X$ is a loop space,
 there is a map $f\:U\to\Omega X$ such that
 $f|_{U_s}=f_s$ for every $s$.
 The maps $\Omega b\circ f$ and $m|_U$ are homotopic
 since they induce the same homomorphism on $\pi_1$.
 \end {demo}

 \begin {claim} {2.9. Lemma.}
 For every $q$,
 there exist
 a cell space $Z$ and
 a map $k\:Z\to U/T$ such that
 $\incl_{(V/T,U/T)}\circ k$ induces an isomorphism
 $\pi_q(Z)\to\pi_q(V/T)$.
 \end {claim}

 \begin {demo} {Proof.}
 For $q\le2$,
 put
 $Z=U/T$,
 $k=\id$.
 We shall construct a single $k$ to serve all $q>2$.

 Take some $s$.
 If $G_s$ if finite,
 let $W_s$ be a space of homotopy type $K(\Z,2)$;
 otherwise,
 let $W_s$ be a point.
 There is a map $a_s\:V_s\to W_s$ with $F(a_s)$ homotopy
 equivalent to $S^1$.
 The map $c_s=a_s|_{T_s}$ is null-homotopic.
 We choose $a_s$ in such a way that $c_s$ is constant.
 Let $j_s\:F(c_s)\to F(a_s)$ be the map induced by
 $\incl_{(V_s,T_s)}$.
 Since $c_s$ is constant,
 $F(c_s)=T_s\times\Omega W_s$.
 Thus $F(c_s)$ is homotopy equivalent to a cell space of
 dimension at most 2.
 By \cite [Lecture~II, Proposition~5] {P},
 $(F(a_s),j_s(F(c_s))$ is a Borsuk pair.
 By Lemma~2.6,
 there is a map $f_s\:F(a_s)\to U_s$ such that
 $\incl_{(V_s,U_s)}\circ f_s$ is rel $j_s(F(c_s))$ homotopic
 to $p(a_s)$.

 $$
 \xymatrix {
 F(c_s)
 \ar[r]^-{j_s}
 \ar[d]^-{p(c_s)} &
 F(a_s)
 \ar[d]^-{p(a_s)}
 \ar[dr]^-{f_s} &
 \\
 T_s
 \ar[r]^-\incl
 \ar[dr]^-{c_s} &
 V_s
 \ar[d]^-{a_s} &
 U_s
 \ar[l]_-\incl \\
 &
 W_s &
 }
 $$

 Put
 $$
 W=W_1\times\ldots\times W_r, \quad
 a=a_1\times\ldots\times a_r\:V\to W, \quad
 c=c_1\times\ldots\times c_r\:T\to W.
 $$
 Since $a|_T$ is constant,
 there is a map $b\:V/T\to W$ such that
 $b\circ\proj_{(V,T)}=a$.
 Let $j\:F(c)\to F(a)$ and $l\:F(a)\to F(b)$ be the maps
 induced by $\incl_{(V,T)}$ and $\proj_{(V,T)}$, respectively.

 We make natural identifications
 $$
 F(a)=F(a_1)\times\ldots\times F(a_r), \quad
 F(c)=F(c_1)\times\ldots\times F(c_r).
 $$
 Then 
 $p(a)=p(a_1)\times\ldots\times p(a_r)$ and
 $j=j_1\times\ldots\times j_r$.
 Put
 $$
 f=f_1\times\ldots\times f_r\:F(a)\to U.
 $$
 The map $\incl_{(V,U)}\circ f$ is rel $j(F(c))$ homotopic to
 $p(a)$.
 Put $g=\proj_{(U,T)}\circ f$.
 We have
 $$
 \incl_{(V/T,U/T)}\circ g=
 \proj_{(V,T)}\circ\incl_{(V,U)}\circ f\sim_{j(F(c))}
 \proj_{(V,T)}\circ p(a)=
 p(b)\circ l.
 $$
 By Lemma~2.7,
 there is a map $h\:F(b)\to U/T$ such that
 $h\circ l=g$ and
 $\incl_{(V/T,U/T)}\circ h\sim p(b)$.
 Since $\pi_q(W)=0$ for $q\ne2$,
 $p(b)$ induces isomorphisms on $\pi_q$, $q>2$.
 Thus the map $\incl_{(V/T,U/T)}\circ h$ does so.
 Thus for any $q>2$ we may
 let $e\:Z\to F(b)$ be a cell approximation and
 put $k=h\circ e$.
 \end {demo}

 $$
 \xymatrix {
 F(c)
 \ar[rrr]^-j
 \ar[dd]^-{p(c)} &
 &
 &
 F(a)
 \ar[rrr]^-l
 \ar[dd]^-{p(a)}
 \ar[dr]_-f
 \ar[drr]^-g &
 &
 &
 F(b)
 \ar[dd]^-{p(b)}
 \ar[dl]_-h \\
 &
 &
 &
 &
 U
 \ar[r]^-\proj
 \ar[dl]^-\incl &
 U/T
 \ar[dr]_-\incl &
 \\
 T
 \ar[rrr]^-\incl
 \ar[drrr]^-c &
 &
 &
 V
 \ar[rrr]^-\proj
 \ar[d]^-a &
 &
 &
 V/T
 \ar[dlll]^-b \\
 &
 &
 &
 W &
 &
 &
 }
 $$


 \subsection * {2.D. Treating a loop space}

 \begin {claim} {2.10. Lemma.}
 Let $X$ be a simply connected space with finitely generated
 homotopy groups.
 Suppose $\pi_2(X)$ has no 2--torsion.
 Then there exist
 a simply connected space $Y$ with finitely generated homotopy
 groups and
 a map $t\:\Omega X\to Y$ inducing split monomorphisms on
 $\pi_q$, $q>1$.
 \end {claim}

 (It follows easily that
 $Y$ may be obtained by attaching 2--cells to $\Omega X$.)

 \begin {demo} {Proof.}
 Put $G=\pi_2(X)$.
 There are
 a space $Q$ of weak homotopy type $K(G,2)$ and
 a map $b\:X\to Q$ inducing an isomorphism on $\pi_2$.
 Consider the piece of the Puppe sequence of $b$:
 $$
 \xymatrix {
 \Omega X
 \ar[r]^-{\Omega b} &
 \Omega Q
 \ar[r]^-j &
 F(b).
 }
 $$
 There is a standard homotopy equivalence
 $e\:\Omega X\to F(j)$ such that
 $p(j)\circ e=\Omega b$.
 Let $V$ be the $K(G,1)$ space considered in 2.C.
 Let $T\subset U\subset V$ be as there.
 There is a weak homotopy equivalence $m\:V\to\Omega Q$.
 By Claim~2.8,
 $m|_U$ lifts (up to homotopy) along $\Omega b$.
 Thus $j\circ m|_U$ is null-homotopic.
 Thus there is a map $h\:C(m|_U)\to F(b)$ such that
 $h\circ i(m|_U)=j$.
 Let $s\:C(m|_T)\to C(m|_U)$ be the map induced by
 $\incl_{(U,T)}$.
 Put $g=h\circ s$.
 Let $r\:F(j)\to F(g)$ be the map induced by $i(m|_T)$.
 Put
 $Y=F(g)$,
 $t=r\circ e$.

 $$
 \xymatrix {
 &
 &
 &
 Y
 \ar@{=}[dr] &
 \\
 &
 \Omega X
 \ar[urr]^-t
 \ar[r]^-e
 \ar[ddr]^-{\Omega b} &
 F(j)
 \ar[rr]^-r
 \ar[dd]^-{p(j)} &
 &
 F(g)
 \ar[dd]^-{p(g)} \\
 &
 U
 \ar[dr]_-{m|_U} &
 &
 &
 \\
 T
 \ar[ur]^-\incl
 \ar[rr]^-{m|_T} &
 &
 \Omega Q
 \ar[rr]^-{i(m|_T)}
 \ar[dr]^-{i(m|_U)}
 \ar[ddr]^-j &
 &
 C(m|_T)
 \ar[dl]_-s
 \ar[ddl]^-g \\
 &
 &
 &
 C(m|_U)
 \ar[d]^-h &
 \\
 &
 &
 &
 F(b) &
 }
 $$

 Let us check the desired properties.
 Since $b$ is 3--connected,
 $F(b)$ is 2--connected.
 Since $m|_T$ is 1--connected,
 $C(m|_T)$ is 1--connected.
 Therefore $g$ is 2--connected.
 Thus $F(g)$ is 1--connected,
 i.~e. $Y$ is simply connected.
 One checks similarly that
 $\pi_q(Y)$ are finitely generated.

 We have the commutative diagram
 $$
 \xymatrix {
 V/T
 \ar[d]_-{\proj_{(V/T,U/T)}} &
 C(\incl_{(V,T)})
 \ar[l]_-\rho
 \ar[r]^-{m'}
 \ar[d]_-\phi &
 C(m|_T)
 \ar[d]_-s \\
 V/U &
 C(\incl_{(V,U)})
 \ar[l]_-\sigma
 \ar[r]^-{m''} &
 C(m|_U),
 }
 $$
 where
 $\phi$ is induced by $\incl_{(U,T)}$,
 $\rho$ and $\sigma$ are the standard homotopy equivalences
 (contractions),
 $m'$ and $m''$ are the weak equivalences induced by $m$.
 Take $q>1$.
 Using Lemma~2.9,
 we get
 a cell space $Z_q$ and
 a map $l_q\:Z_q\to V/T$ inducing an isomorphism on $\pi_q$ and
 such that
 $\proj_{(V/T,U/T)}\circ l_q$ is constant.
 Using the diagram,
 we get a map $v_q\:Z_q\to C(m|_T)$ inducing an isomorphism on
 $\pi_q$ and such that
 $s\circ v_q$ is null-homotopic.
 Since $g=h\circ s$,
 $g\circ v_q$ is null-homotopic.
 Thus $v_q$ lifts along $p(g)$:
 there is a map $w_q\:Z_q\to F(g)$ such that
 $p(g)\circ w_q=v_q$.
 Consider the commutative diagram
 $$
 \xymatrix {
 &
 0
 \ar@{=}[d] &
 &
 &
 0
 \ar@{=}[d] \\
 &
 \pi_{q+1}(\Omega Q)
 \ar[r]^-{j_*}
 \ar[d]^-{i(m|_T)_*} &
 \pi_{q+1}(F(b))
 \ar[r]^-{\d_j}
 \ar@{=}[d] &
 \pi_q(F(j))
 \ar[r]^-{p(j)_*}
 \ar[d]^-{r_*} &
 \pi_q(\Omega Q)
 \ar[d]^-{i(m|_T)_*} \\
 \pi_{q+1}(F(g))
 \ar[r] &
 \pi_{q+1}(C(m|_T))
 \ar[r]^-{g_*} &
 \pi_{q+1}(F(b))
 \ar[r]^-{\d_g} &
 \pi_q(F(g))
 \ar[r] &
 \pi_q(C(m|_T)) \\
 &
 \pi_{q+1}(Z_{q+1})
 \ar[ul]^-{w_{q+1\,*}}
 \ar[u]_-{v_{q+1\,*}} &
 &
 &
 \pi_q(Z_q),
 \ar[ul]^-{w_{q\,*}}
 \ar[u]_-{v_{q\,*}}
 }
 $$
 where the rows are pieces of the exact sequences of the maps
 $j$ and $g$.
 We see that
 $\d_j$ is an isomorphism and
 $\d_g$ is a split monomorphism.
 Thus $r_*$ is a split monomorphism.
 Since $t=r\circ e$, where $e$ is a homotopy equivalence,
 $t$ induces a split monomorphism on $\pi_q$.
 \end {demo}

 \begin {claim} {2.11. Corollary.}
 Let $X$ be a simply connected space with finitely generated
 homotopy groups.
 Then there exist
 a simply connected space $Y$ with finitely generated homotopy
 groups and
 a map $t\:\Omega X\to Y$ such that
 for every $q>1$
 the homomorphism $t_*\otimes\id\:
 \pi_q(\Omega X)\otimes\Z[1/2]\to\pi_q(Y)\otimes\Z[1/2]$ is a
 split monomorphism.
 \end {claim}

 \begin {demo} {Proof.}
 By attaching 3-- and 4--cells to $X$,
 we may quotient $H_2$ by 2--torsion
 without changing $H_q$ for $q>2$.
 By Serre,
 $\pi_q\otimes\Z[1/2]$ are not changed.
 Then apply Lemma~2.10.
 \end {demo}


 \section * {3. Dimension descent}


 \subsection * {3.A. Concentrating spheroids
                near a codimension 2 net in $I^n$} 

 \begin {claim} {3.1. Lemma.}
 Let
 $X$ be a simply connected space,
 $\Gamma$ be a finite open cover of $I^n$ ($n>1$).
 Then there exist
 a finite open cover $\Delta$ of $I^{n-1}$,
 a function $\lambda\:\Delta\to\Gamma(n-1)$, and
 a function $\Phi\:\Pi_n(X)\to\Pi_n(X)$ such that
 (1) $[\Phi(a)]=[a]$ for all $a\in\Pi_n(X)$ and
 (2) for $F\in\Delta$, $a,a'\in\Pi_n(X)$,
 the implication holds
 $$
 a|_{\lambda(F)}=a'|_{\lambda(F)} \quad
 \Rightarrow \quad
 \Phi(a)|_{F\times I}=\Phi(a')|_{F\times I}.
 $$
 \end {claim}

 \begin {demo} {Proof.}
 Choose $\epsilon>0$ such that
 every open ball $B(u,3\epsilon)$ ($u\in I^n$) is contained in
 some $E\in\Gamma$.
 Choose a rectilinear triangulation $T$ of $I^n$ with simplices
 of diameter less than $\epsilon$.
 Let $T'$ be the barycentric subdivision of $T$.
 By $\cosk^2T$ we denote the union of all simplices of $T'$
 that do not intersect $\sk_1T$.
 (So $\cosk^2T$ is a polyhedron of dimension $n-2$;
 its complement collapses to $\sk_1T$.)

 There is a homeomorphism $l\:I^n\to I^n$ preserving simplices
 of $T$ and such that
 $$
 |v\times I\cap l^{-1}(\cosk^2T)|\le n-1
 $$
 for every $v\in I^{n-1}$.
 (Such $l$ is obtained by generic simplex-wise perturbation of
 the identity map.)
 The map $l$ preserves $\d I^n$ and has degree 1.

 For $v\in I^{n-1}$,
 let $2\rho(v)$ be the least positive distance between
 $v\times I$ and $l^{-1}(t)$, where
 $t$ runs over simplices of $T'$.
 Let $\Delta$ be a finite cover of $I^{n-1}$ formed by some of
 the balls $B(v,\rho(v))$, $v\in I^{n-1}$.

 Let us construct the function $\lambda$.
 Take $F\in\Delta$.
 We have $F=B(v,\rho(v))$ for some $v\in I^{n-1}$.
 For each
 $$
 u\in v\times I\cap l^{-1}(\cosk^2T)
 $$
 (there are at most $n-1$ such $u$),
 choose $E\in\Gamma$ containing $B(u,3\epsilon)$.
 Put $\lambda(F)$ be the union of these $E$.

 Choose $\delta>0$ such that
 $\delta\le\epsilon$ and
 $\delta\le\rho(v)$ for $B(v,\rho(v))\in\Delta$.
 Let $Q$ be the open $\delta$--neighbourhood of
 $l^{-1}(\cosk^2T)$.
 Put $P=l(I^n\setminus Q)$.
 The set $P$ is closed and does not intersect $\cosk^2T$.
 Thus there is a map $k\:I^n\to I^n$ preserving the simplices
 of $T$ with $k(P)\subset\sk_1T$.
 The map $k$ preserves $\d I^n$ and has degree 1.

 Take $a\in\Pi_n(X)$.
 We homotop $a$ to get a map $\Theta(a)$ taking $\sk_1T$ to the
 basepoint.
 The homotopy is constructed by induction on the skeleta of $T$.
 There are no obstructions
 since $X$ is simply connected.
 To extend the homotopy to some simplex $s$ of $T$,
 we need to know only $a|_s$ and the homotopy constructed on
 $\d s$
 (we need no information from the outside of $s$).

 Therefore there is a function $\Theta\:\Pi_n(X)\to\Pi_n(X)$
 such that
 (a) $[\Theta(a)]=[a]$ for all $a\in\Pi_n(X)$,
 (b) for every simplex $s$ of $T$ and any $a,a'\in\Pi_n(X)$,
 the implication holds
 $$
 a|_s=a'|_s \quad
 \Rightarrow \quad
 \Theta(a)|_s=\Theta(a')|_s,
 $$
 and 
 (c) for any $a\in\Pi_n(X)$,
 $\Theta(a)$ takes $\sk_1T$ to the basepoint.

 For $a\in\Pi_n(X)$,
 put $\Phi(a)=\Theta(a)\circ k\circ l$.
 The property (1) is obvious.
 Let us check the property (2).
 Take
 $F\in\Delta$ and
 $a,a'\in\Pi_n(X)$ such that
 $a|_{\lambda(F)}=a'|_{\lambda(F)}$.
 We should show that
 $\Phi(a)|_{F\times I}=\Phi(a')|_{F\times I}$.
 Take $u_0\in F\times I$.
 Let us show that $\Phi(a)(u_0)=\Phi(a')(u_0)$.

 If $l(u_0)\in P$, then
 $k(l(u_0))\in\sk_1T$ and
 thus $\Phi(a)(u_0)=\Phi(a')(u_0)=x_0$,
 where $x_0\in X$ is the basepoint.
 Consider the converse case.
 Since $P=l(I^n\setminus Q)$,
 we have $u_0\in Q$.
 Thus $\dist(u_0,l^{-1}(\cosk^2T))<\delta$.
 Thus there is a simplex $t$ of $\cosk^2T$ such that
 $\dist(u_0,l^{-1}(t))<\delta$.
 We have $F=B(v,\rho(v))$ for some $v\in I^{n-1}$.
 We have
 $$
 \dist(v\times I,l^{-1}(t))\le
 \dist(v\times I,u_0)+\dist(u_0,l^{-1}(t))<
 \rho(v)+\delta\le
 2\rho(v).
 $$
 By definition of $\rho(v)$, this means that
 $\dist(v\times I,l^{-1}(t))=0$.
 Thus there is a point $u\in v\times I\cap l^{-1}(t)$.
 Since $l$ preserves simplices of $T$,
 $\diam l^{-1}(t)<\epsilon$.
 We have
 $$
 \dist(u_0,u)\le
 \dist(u_0,l^{-1}(t))+\diam l^{-1}(t)<
 \delta+\epsilon\le
 2\epsilon.
 $$
 The point $u_0$ belongs to some simplex $s$ of $T$.
 We have
 $$
 s\subset
 B(u_0,\epsilon)\subset
 B(u,3\epsilon)\subset
 \lambda(F).
 $$
 Thus $a|_s=a'|_s$.
 Thus $\Theta(a)|_s=\Theta(a')|_s$.
 Since $k$ and $l$ preserve $s$,
 $\Phi(a)|_s=\Phi(a')|_s$,
 which suffices.
 \end {demo}


 \subsection * {3.B. Functionals of finite degree}

 For a space $X$,
 the formula
 $\Xi(a)(t_1,\ldots,t_{n-1})(t)=a(t_1,\ldots,t_{n-1},t)$
 ($t_1,\ldots,t_{n-1},t\in I$)
 defines a bijection $\Xi\:\Pi_n(X)\to\Pi_{n-1}(\Omega X)$,
 which we call \textit {standard}.
 The induced isomorphism $\xi\:\pi_n(X)\to\pi_{n-1}(\Omega X)$
 we also call standard.
 For an open set $F\subset I^{n-1}$, 
 the bijection $\Xi_F\:\Pi(F\times I,X)\to\Pi(F,\Omega X)$
 defined by that formula is also called standard.

 \begin {claim} {3.2. Corollary.}
 Let
 $X$ be a simply connected space,
 $L$ be an abelian group,
 $g\:\Pi_{n-1}(\Omega X)\to L$ be a homotopy invariant
 functional ($n>1$).
 Let $\Xi\:\Pi_n(X)\to\Pi_{n-1}(\Omega X)$ be the standard
 bijection.
 Then $\deg g\circ\Xi\le(n-1)\deg g$.
 \end {claim}

 \begin {demo} {Proof.}
 Suppose $\deg g\le r$.
 Let us show that $\deg g\circ\Xi\le(n-1)r$.
 Let $\Gamma$ be a finite open cover of $I^n$.
 By Lemma~3.1,
 there are
 a finite open cover $\Delta$ of $I^{n-1}$,
 a function $\lambda\:\Delta\to\Gamma(n-1)$, and
 a function $\Phi\:\Pi_n(X)\to\Pi_n(X)$ satisfying the conditions
 (1), (2) of the lemma.

 For each $F\in\Delta(r)$,
 choose a decomposition $F=F_1\cup\ldots\cup F_s$,
 $0\le s\le r$, with $F_1,\ldots,F_s\in\Delta$ and
 put $\mu(F)=\lambda(F_1)\cup\ldots\cup\lambda(F_s)$.
 So we have a function $\mu\:\Delta(r)\to\Gamma((n-1)r)$.
 It follows from the condition (2) that
 for $F\in\Delta(r)$, $a,a'\in\Pi_n(X)$
 the implication holds
 $$
 a|_{\mu(F)}=a'|_{\mu(F)} \quad
 \Rightarrow \quad
 \Phi(a)|_{F\times I}=\Phi(a')|_{F\times I}.
 $$
 Thus for every $F\in\Delta(r)$
 there is a function $\Phi_F\:\Pi(\mu(F),X)\to\Pi(F\times I,X)$
 such that
 $\Phi_F(a|_{\mu(F)})=\Phi(a)|_{F\times I}$
 for all $a\in\Pi_n(X)$.
 For every $F\in\Delta(r)$,
 we have the commutative diagram
 $$
 \xymatrix {
 \Pi_n(X)
 \ar[r]^-\Phi
 \ar[d] &
 \Pi_n(X)
 \ar[r]^-\Xi
 \ar[d] &
 \Pi_{n-1}(\Omega X)
 \ar[d] \\
 \Pi(\mu(F),X)
 \ar[r]^-{\Phi_F} &
 \Pi(F\times I,X)
 \ar[r]^-{\Xi_F} &
 \Pi(F,\Omega X),
 }
 $$
 where
 $\Xi_F$ is the standard bijection and
 the vertical arrows are the restriction functions.

 Since $\deg g\le r$,
 there are functionals $g_F\:\Pi(F,\Omega X)\to L$,
 $F\in\Delta(r)$, such that
 $$
 g(b)=\sum_{F\in\Delta(r)}g_F(b|_F)
 $$
 for all $b\in\Pi_{n-1}(\Omega X)$.
 For $E\in\Gamma((n-1)r)$,
 define a functional $f_E\:\Pi(E,X)\to L$ by
 $$
 f_E(a)=\sum_{F\in\mu^{-1}(E)}g_F(\Xi_F(\Phi_F(a))).
 $$
 For $a\in\Pi_n(X)$,
 we have
 \begin {multline*}
 g(\Xi(a))=
 g(\Xi(\Phi(a)))=
 \sum_{F\in\Delta(r)}g_F(\Xi(\Phi(a))|_F)=
 \sum_{F\in\Delta(r)}g_F(\Xi_F(\Phi_F(a|_{\mu(F)})))= \\
 =\sum_{E\in\Gamma((n-1)r)}
 \sum_{F\in\mu^{-1}(E)}g_F(\Xi_F(\Phi_F(a|_E)))=
 \sum_{E\in\Gamma((n-1)r)}f_E(a|_E).
 \end {multline*}
 \end {demo}

 \begin {claim} {3.3. Lemma.}
 Let
 $X$ and $Y$ be spaces,
 $t\:X\to Y$ be a map.
 Let
 $L$ be an abelian group,
 $g\:\Pi_n(Y)\to L$ be a functional.
 Then $\deg g\circ t_\#\le\deg g$.
 \end {claim}

 \begin {demo} {}
 (Obvious.)
 \end {demo}

 \begin {claim} {3.4. Lemma.}
 Let $X$ be a space.
 Define $l\:\Pi_n(X)\to H_n(X)$ by $l(a)=a_*(u)$,
 where $u\in H_n(I^n,\d I^n)$ is the fundamental class.
 Then $\deg l\le1$.
 \end {claim}

 \begin {demo} {Proof.}
 Let $\Gamma$ be a finite open cover of $I^n$.
 Represent $u$ by a (singular) cycle $U\in Z_n(I^n,\d I^n)$
 subordinate to $\Gamma$:
 $$
 U=\sum_{E\in\Gamma}\incl_{(I^n,E)\,\#}(U_E),
 $$
 where $U_E\in C_n(E,E\cap\d I^n)$ are some chains.

 The subgroup $Z_n(X)\subset C_n(X)$ is a direct summand.
 Thus there is a homomorphism $k\:C_n(X)\to H_n(X)$ such that
 $k(T)=[T]$
 for all $T\in Z_n(X)$.
 For $E\in\Gamma$,
 define $l_E\:\Pi(E,X)\to H_n(X)$ by $l_E(a)=k(a_\#(U_E))$.
 For $a\in\Pi_n(X)$,
 we have
 $$
 l(a)=
 [a_\#(U)]=
 k(a_\#(U))=
 \sum_{E\in\Gamma}k((a|_E)_\#(U_E))=
 \sum_{E\in\Gamma}l_E(a|_E).
 $$
 \end {demo}

 \begin {demo} {\bf Proof of Theorem~1.1.}
 Induction on $n$.
 For $n=2$,
 consider the commutative diagram
 $$
 \xymatrix {
 &
 \Pi_2(X)
 \ar[dl]_-l
 \ar[d]^-p
 \ar[dr]^-q \\
 H_2(X) &
 \pi_2(X)
 \ar[l]_-h
 \ar[r]^-m &
 \pi_2(X)\otimes\Z[1/2],
 }
 $$
 where
 $p$ is the natural projection,
 $m$ is defined by $m(v)=v\otimes1$,
 $h$ is the Hurewicz isomorphism, and
 $l$ is as in Lemma~3.4,
 thus $\deg l\le1$.
 We see that $\deg p\le1$ and
 thus $\deg q\le1$.

 Take $n>2$.
 By Corollary~2.11,
 there are
 a simply connected space $Y$ with finitely generated homotopy
 groups and
 a map $t\:\Omega X\to Y$ such that
 for every $q>1$
 the homomorphism $t_*\otimes\id\:
 \pi_q(\Omega X)\otimes\Z[1/2]\to\pi_q(Y)\otimes\Z[1/2]$ is a
 split monomorphism.
 Consider the commutative diagram
 $$
 \xymatrix {
 \Pi_n(X)
 \ar[r]^-\Xi
 \ar[d]^-q &
 \Pi_{n-1}(\Omega X)
 \ar[r]^-{t_\#}
 \ar[d]^-{q'} &
 \Pi_{n-1}(Y)
 \ar[d]^-{q''} \\
 \pi_n(X)\otimes\Z[1/2]
 \ar[r]^-{\xi\otimes\id} &
 \pi_{n-1}(\Omega X)\otimes\Z[1/2]
 \ar[r]^-{t_*\otimes\id} &
 \pi_{n-1}(Y)\otimes\Z[1/2],
 }
 $$
 where
 $\Xi$ is the standard bijection,
 $\xi$ is the standard isomorphism,
 $q'$ and $q''$ are defined similarly to $q$.
 Since $\xi\otimes\id$ is an isomorphism,
 $\deg q=\deg(\xi\otimes\id)\circ q=\deg q'\circ\Xi$.
 By Corollary~3.2,
 $\deg q'\circ\Xi\le(n-1)\deg q'$.
 Since $t_*\otimes\id$ is a split monomorphism,
 $\deg q'=\deg(t_*\otimes\id)\circ q'=\deg q''\circ t_\#$.
 By Lemma~3.3,
 $\deg q''\circ t_\#\le\deg q''$.
 By induction hypothesis,
 $\deg q''\le(n-2)!$.
 Therefore $\deg q\le(n-1)!$.
 \end {demo}
 


 {\noindent \tt ssp@pdmi.ras.ru}

 {\noindent \tt http://www.pdmi.ras.ru/\~{}ssp}

 \end {document}